\documentclass[10pt,a4paper]{amsart}

\author[M.~Truter]{Marc~Truter {\large\orcidlink{0009-0009-5209-2238}}}
\address{Mathematics Institute\\University of Warwick\\Coventry\\CV4 7AL\\UK}
\email{Marc.Truter@warwick.ac.uk}
\date{}

\usepackage[margin=2.5cm]{geometry}
\usepackage{amsmath}
\usepackage{amsfonts} 
\usepackage{mathrsfs}
\usepackage{amssymb}
\usepackage{graphicx}
\usepackage[font={small,it}]{caption}
\usepackage[font={small,it}]{subcaption}

\usepackage{enumitem}
\usepackage[hidelinks]{hyperref}
\usepackage{orcidlink}
\hypersetup{hidelinks}
\usepackage{booktabs}
\usepackage{colortbl}

\usepackage[numbers,sort&compress]{natbib}

\begin{document}

\title[Terminal Fano 4-fold Hypersurfaces
via Deep Reinforcement Learning]{Terminal Fano 4-fold Hypersurfaces \\ via Deep Reinforcement Learning}

\begin{abstract}
{\small We construct $46,992$ new families of Fano $4$-fold hypersurfaces $X_d\subset\mathbb(a_1,\ldots,a_6)$ of Fano index~$1$ with $d>200$, whose general members have terminal singularities. Unlike the $11{,}617$ quasismooth families classified by Brown and Kasprzyk, these are nonquasismooth, for which determining terminality can require the analysis of hyperquotient singularities. An exhaustive search has classified all nonquasismooth families with $d\leq200$, but extending this classification to higher degrees becomes computationally prohibitive. We therefore exploit the sparsity and spatial clustering of terminal examples using two best-first searches over $\mathbb{Z}^6$: one with a fixed heuristic, and the other with a neural network trained online by deep reinforcement learning. Together, the searches reach degrees up to $3{,}234$. The neural network additionally discovers an average of $79$ hard-to-reach families per run. For the most isolated examples, our conservative distance-based model gives an estimated lower bound of nearly $40\%$ of the fixed search's total step budget to discover a single such family. An ablation study demonstrates that this advantage arises from the neural network rather than stochastic exploration alone. The complete list of families is publicly available.}
\end{abstract}


\maketitle

\section{Introduction}

\begin{figure}[!h]
    \centering
    \includegraphics[width=0.36\linewidth]{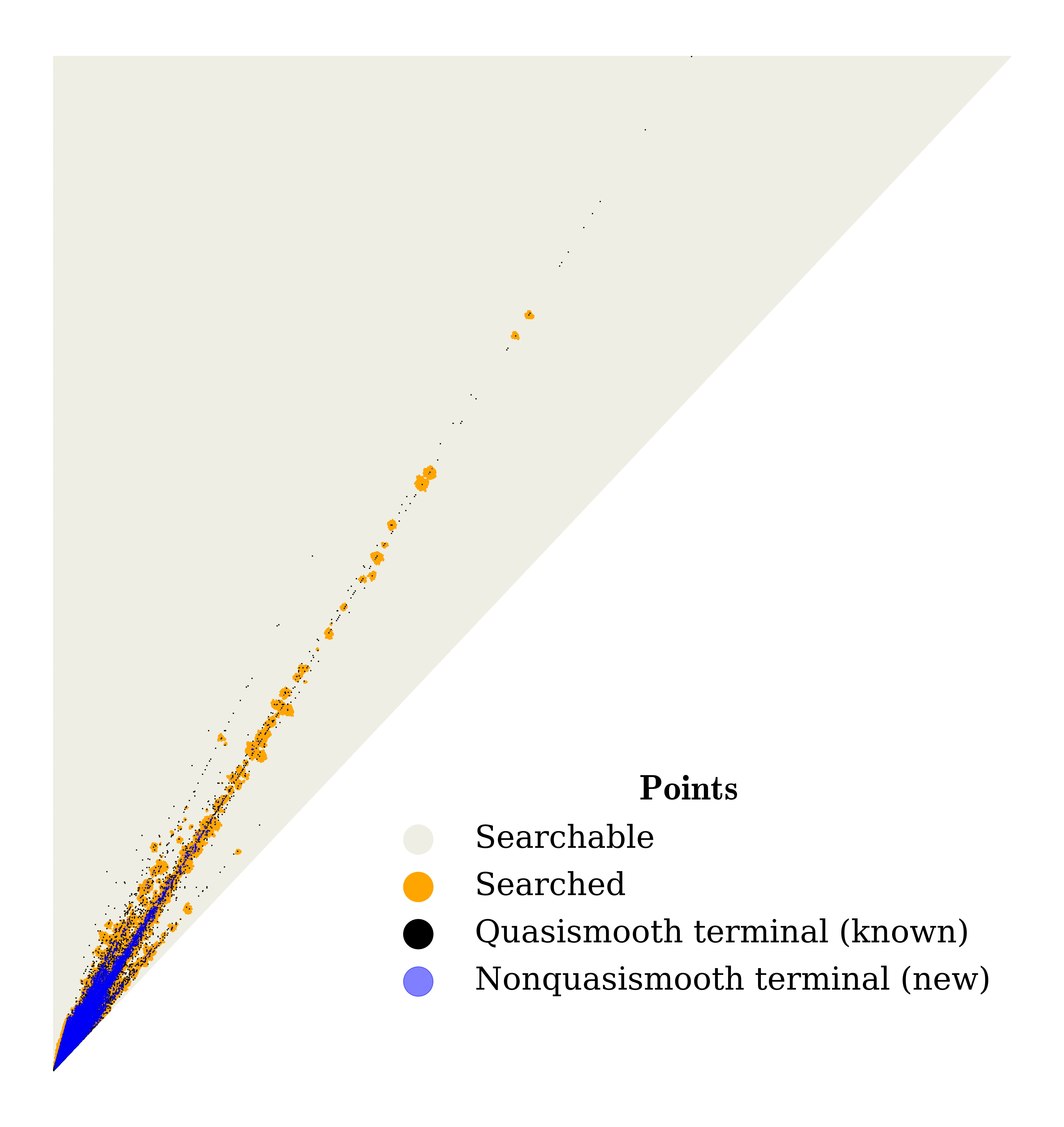}
    \caption{The dynamic heuristic search over $\mathbb{Z}^6$,
    projected onto two dimensions. Alongside the classified
    quasismooth families, the search discovers previously
    unknown nonquasismooth terminal examples.
    See Figure~\ref{fig:search_dynamic} for details.}
    \label{fig:search_dynamic_hero}
\end{figure}

We use heuristic search algorithms to discover new families of nonquasismooth terminal Fano $4$-fold hypersurfaces of index~$1$ with degree $d>200$. Figure~\ref{fig:search_dynamic_hero} illustrates how a neural network guided search explores the space of candidate weight systems and discovers previously unknown families.

\subsection{Main results}
\label{sec:main_results}

The minimal model program provides a framework for the birational classification of algebraic varieties, in which Fano varieties, Calabi--Yau varieties, and varieties of general type play central roles. It allows terminal singularities, a particularly mild class of singularities arising naturally in higher-dimensional birational geometry \citep{YPG}\citep[\S2.2]{Reid2002}. Birkar~\citep{B} proved that Fano varieties with terminal singularities form finitely many families in each fixed dimension, making their complete classification, a `periodic table', a finite problem. Such classifications are well understood in dimensions $1$ and $2$, while substantial progress has been made in dimension $3$ \citep{BKR,I1,I2,MM}. Notably, Reid produced the list of $95$ Fano $3$-fold hypersurface families of index $1$ \citep[\S16.6]{F}, and Iano-Fletcher found the corresponding list of $85$ codimension-$2$ weighted complete intersection families \citep[\S16.7]{F}. In dimension~$4$, however, the classification remains in its infancy, with the Fano index~$1$ case wide open even for smooth varieties~\citep{CKP}.

We consider families of Fano $4$-fold hypersurfaces $X_d\subset\mathbb{P}(a_1,\ldots,a_6)$ of index $1$, whose general members have terminal singularities. By adjunction, their degree satisfies $d=\left(\sum a_i\right)-1$, so candidate weight systems are encoded by integer vectors $(a_1,\ldots,a_6)\in\mathbb{Z}^6$, where $1\leq a_1\leq\cdots\leq a_6$. Brown and Kasprzyk~\citep{BK} classified the $11{,}617$ quasismooth families, recorded in the Graded Ring Database~\citep{GRDB}. Quasismoothness considerably simplifies the problem: singularities are inherited from the ambient weighted projective space, making terminality comparatively inexpensive to test, while strong restrictions on the weights permit a terminating tree-search algorithm. An exhaustive search classifies all nonquasismooth terminal families of index~$1$ up to degree~$200$~\citep{MTFano, Thesis}, with the corresponding classification data and source code
archived on Zenodo~\citep{MTFano200, MTFanoExhaustiveSoftware}.

\begin{table}[!h]
    \centering
    \small
    \setlength{\arrayrulewidth}{3pt}
    \arrayrulecolor{white}
    \begin{tabular}{@{}lllll@{}}
        \rowcolor{gray!30}
        \textbf{Method}
        & \textbf{Families}
        & $\boldsymbol{d>200}$
        & \textbf{Max degree}
        & \textbf{Hard-to-reach} \\ \hline
        \rowcolor{gray!10}
        Exhaustive ($d\leq200$)
            & $77,353$ & $0$ & $200$ & --- \\ \hline
        \rowcolor{gray!20}
        Fixed heuristic
            & $113{,}996$ & $41{,}376$
            & $3{,}234$ & $0$ \\
        \rowcolor{gray!20}
        Dynamic heuristic
            & $81{,}823\pm5{,}831$ & $16{,}632\pm2{,}685$
            & $2{,}659\pm417$ & $79\pm9$ \\ 
        \rowcolor{gray!30}
        \textbf{Combined union}
            & $\boldsymbol{121{,}980}$ & $\boldsymbol{46{,}992}$
            & $\boldsymbol{3{,}234}$ & --- \\
    \end{tabular}
    \caption{Summary of heuristic search results compared with exhaustive classification for $d\leq200$. Counts refer to distinct nonquasismooth terminal families. Dynamic results are mean $\pm$ std over five runs, while the union consists of all distinct families found across the heuristic searches. ``Hard-to-reach'' denotes families with $d>200$ and $D(p)\geq10$, as defined in \S\ref{sec:search_results}.}
    \label{tab:summary}
\end{table}

In the present work, we extend the known list using fixed and dynamic heuristic searches, yielding $46{,}992$ previously unknown nonquasismooth Fano $4$-fold hypersurface families of index~$1$ with $d>200$, each certified to have a terminal general member. Their degrees reach $3{,}234$. Across all degrees, the searches find $121{,}980$ distinct nonquasismooth terminal families. Table~\ref{tab:summary} summarises these results. Although the fixed heuristic discovers more families overall, the dynamic heuristic reaches more isolated regions of weight space, finding on average $79\pm9$ hard-to-reach families per run. For the most isolated examples, our conservative distance-based model estimates lower bounds approaching $40\%$ of the fixed search's step budget to discover a single family, as detailed in \S\ref{sec:search_results}. An ablation study confirms that this advantage arises from learning rather than stochastic exploration alone, as discussed in \S\ref{sec:ablation}.

Selected examples of the newly discovered higher-degree families are presented in Table~\ref{tab:search_results}. The complete list of $121{,}980$ families, together with their weights, singularity baskets and anticanonical volumes, is archived on Zenodo~\citep{MTFanoSporadic}, alongside the source code~\citep{MTFanoHeuristicSoftware}. As an illustration, the database records for $X_{136}\subset\mathbb{P}(1,1,5,25,37,68)$ two singular points $\frac{1}{37}(1,1,5,31)$ and $\frac{1}{25}(1,1,5,12,18;11)$, a singular curve $\frac{1}{5}(0,1,2,3)$, and anticanonical volume $(-K_X)^4=\frac{2}{4625}$.

\subsection{The computational obstruction}
\label{sec:computational_obstruction}

As noted in \S\ref{sec:main_results}, an exhaustive search conducted in separate work classifies all nonquasismooth terminal Fano $4$-fold hypersurface families of index~$1$ up to degree~$200$, yielding $77{,}353$ families. Extending this classification to higher degrees poses substantial computational challenges. The number of candidate weight vectors of degree $d$ is $O(d^5)$, giving $O(D^6)$ candidates over all degrees $d\leq D$. Moreover, unlike in the quasismooth case, determining terminality can require costly analysis of hyperquotient singularities. Together, these difficulties make exhaustive enumeration at higher degrees computationally prohibitive.

\begin{figure}[!h]
    \centering
    \includegraphics[width=0.51\linewidth]{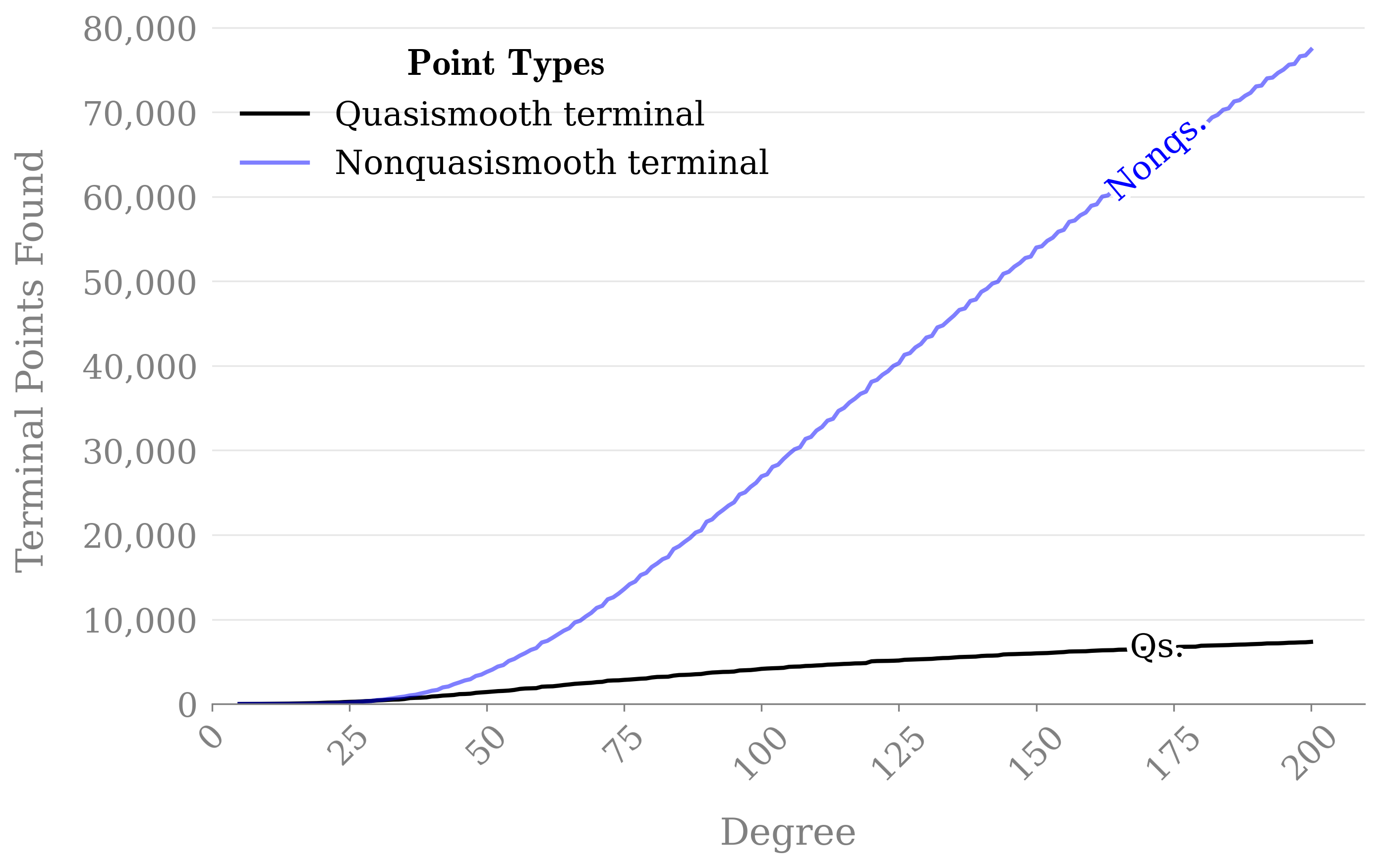}
    \caption{Cumulative numbers of quasismooth and nonquasismooth
    terminal Fano $4$-fold hypersurface families of degree at
    most $d$, found by exhaustive search.}
    \label{fig:qs_vs_nonqs}
\end{figure}

\subsection{The algorithmic contribution}

\begin{figure}[!h]
    \centering
    \includegraphics[width=0.8\linewidth]{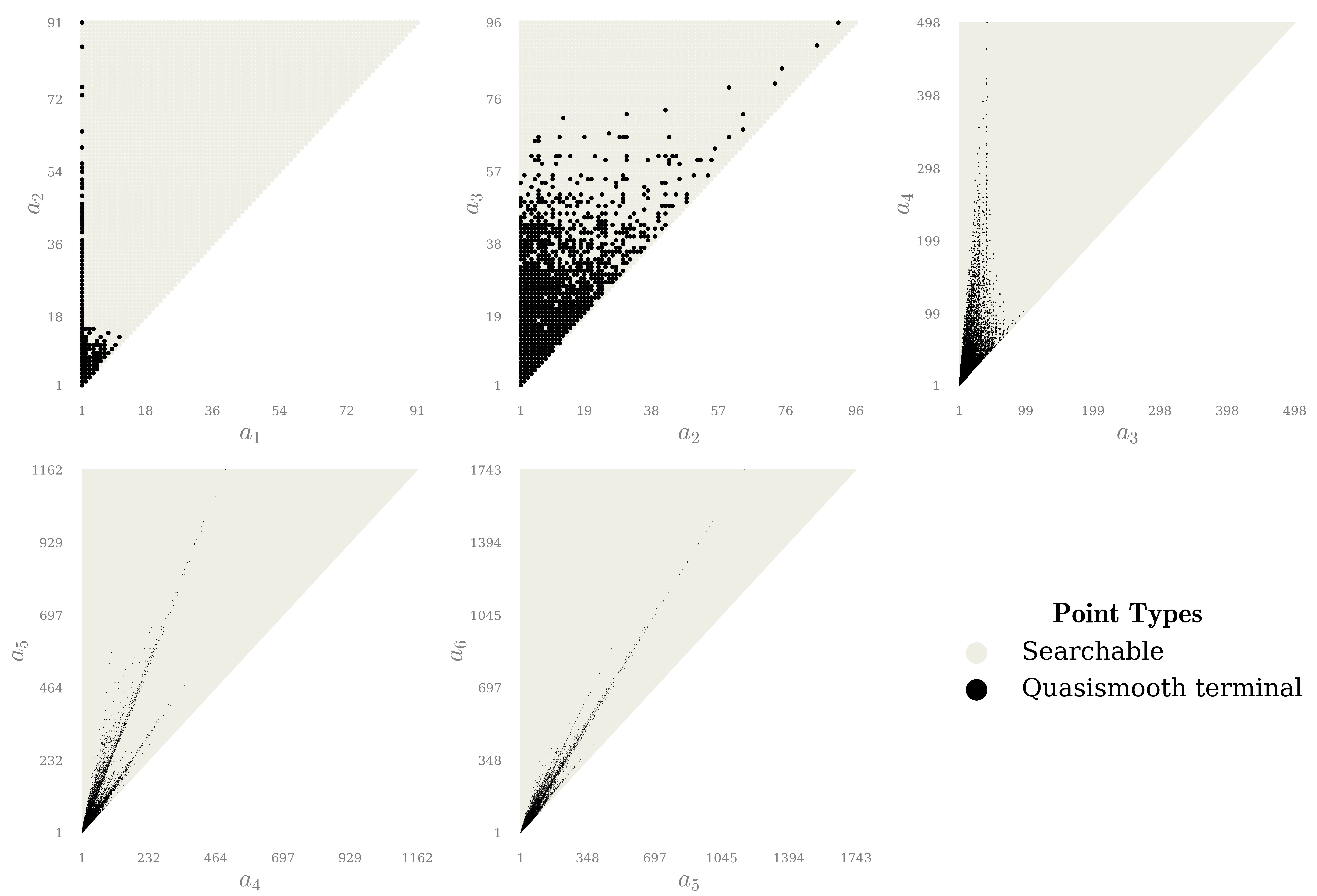}
    \caption{The $11{,}617$ quasismooth terminal Fano $4$-fold
    hypersurface families represented by weight vectors in
    $\mathbb{Z}^6$, projected onto consecutive coordinate
    pairs $(a_i,a_{i+1})$. The projections illustrate their
    concentration in weight space.}
    \label{fig:class_qs}
\end{figure}

The quasismooth classification exhibits a striking structure in weight space: the projections in Figure~\ref{fig:class_qs} suggest that terminal examples are spatially clustered. This motivates searching from previously discovered examples.

We recast the classification as a search over $\mathbb{Z}^6$, where terminal weight systems are \emph{reward points}. We exploit this structure using two best-first search algorithms, defined in \S\ref{sec:integer_lattice_search}. The \emph{fixed heuristic} is deterministic, prioritising candidates by proximity to previously discovered rewards. The \emph{dynamic heuristic} instead uses a compact neural network, trained online by deep reinforcement learning via temporal difference learning, together with stochastic exploration. The network smooths the high-variance reward signal, guiding the search towards regions the fixed heuristic struggles to reach. In both methods, the heuristics select candidate weight vectors, while terminality is tested independently using an algorithmic implementation of the relevant mathematical criteria.

This work builds on research applying data science and machine learning to mathematical data, particularly in algebraic geometry \citep{BYHHJL,CKVI,CKVII}, with subsequent applications across mathematics \citep{MAI,H2024}.

\subsection{Paper outline}

\S~\ref{sec:background} recalls the mathematical background and terminality criteria. \S~\ref{sec:integer_lattice_search} introduces the search algorithms and their implementation. \S~\ref{sec:results} presents the computational results, including a distance based analysis of hard-to-reach examples, and an ablation study.

\section{Mathematical background}
\label{sec:background}

To ground the general construction, we first illustrate it with a classical example, elliptic curves.  The family of all elliptic curves is given by $X_6\subset\mathbb{P}(1,2,3)$. The ambient \textit{weighted projective space}, $\mathbb{P}(1,2,3)=(\mathbb{C}^3\backslash\{0\}) / \mathbb{C}^*$, where $\lambda\in\mathbb{C}^*$ acts on $\mathbb{C}^3\backslash\{0\}$ with coordinates $(x,y,z)$ via $\lambda\cdot(x,y,z)=(\lambda x,\lambda^2 y,\lambda^3 z)$. A curve in the family $X_6:(f_6=0)$ is the set of solutions of a homogeneous polynomial $f_6$ of degree $6$ which must be of the form
\[
f_6=c_1z^2+c_2y^3+c_3x^6+c_4x^4y+c_5x^2y^2+c_6x^3z+c_7xyz
\]
for some $c_1,\ldots,c_7\in \mathbb{C}$, noting that $x$ has weight $1$, $y$ has weight $2$, and $z$ has weight $3$, so that each term does indeed have weight~6. Therefore, the family $X_6$ given by all possible equations $f_6$ is parametrised by its coefficients $[c_1:\cdots:c_7]\in \mathbb{P}^6$. 

Extending the same construction to any weight $d$ and dimension $n$, we can define families of $n$-dimensional hypersurfaces 
\[
X_d\colon (f_d=0)\subset\mathbb{P}(a_1,\ldots,a_{n+2})
\]
for weights $1\le a_1\le\hdots\le a_{n+2}$. As with the elliptic curve example, the family is parametrised by $\mathbb{P}^{N-1}$, where $N$ is the number of monomials of degree $d$ in weights $a_1,\ldots,a_{n+2}$. The \textit{adjunction number} is defined as
\[
\alpha=d-\sum^{n+2}_{i=1} a_i 
\]
and $X_d$ is Fano if $\alpha<0$, Calabi-Yau if $\alpha=0$, and general type if $\alpha>0$. For example, the elliptic curves $X_6\subset\mathbb{P}(1,2,3)$ have $\alpha=6-(1+2+3)=0$, and they are Calabi-Yau varieties. When $\alpha<0$, we define the \textit{Fano index} to be $i_X=-\alpha>0$. 

In this paper, we consider the main case of terminal Fano $4$-fold hypersurfaces, those of Fano index $i_X=1$. By fixing the Fano index, we have $d=\left(\sum a_i\right)-1$, and so may encode the data as an integer vector $(a_1,\ldots,a_6)\in\mathbb{Z}^6$ bounded by $1\le a_1\le\hdots\le a_{n+2}$. 

\begin{figure}[!h]
    \centering
    \begin{subfigure}[b]{0.12\linewidth}
        \centering
        \includegraphics[width=\linewidth]{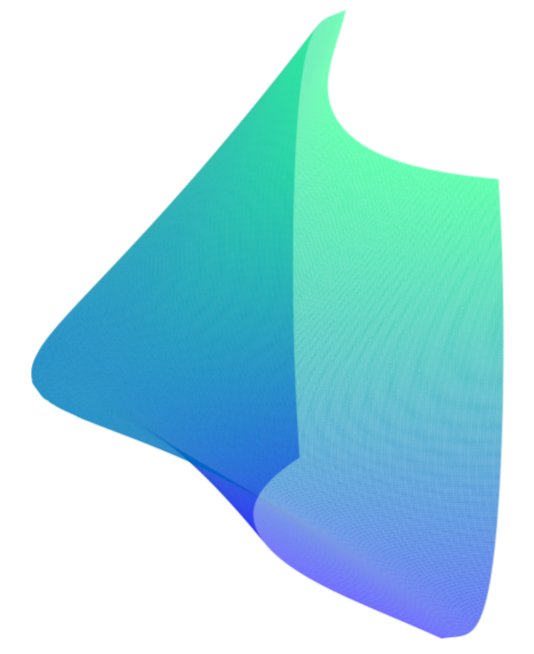}
        \caption{}
        \label{fig:singularities_qs}
    \end{subfigure}
    \hspace{0.15\linewidth}
    \begin{subfigure}[b]{0.2\linewidth}
        \centering
        \includegraphics[width=\linewidth]{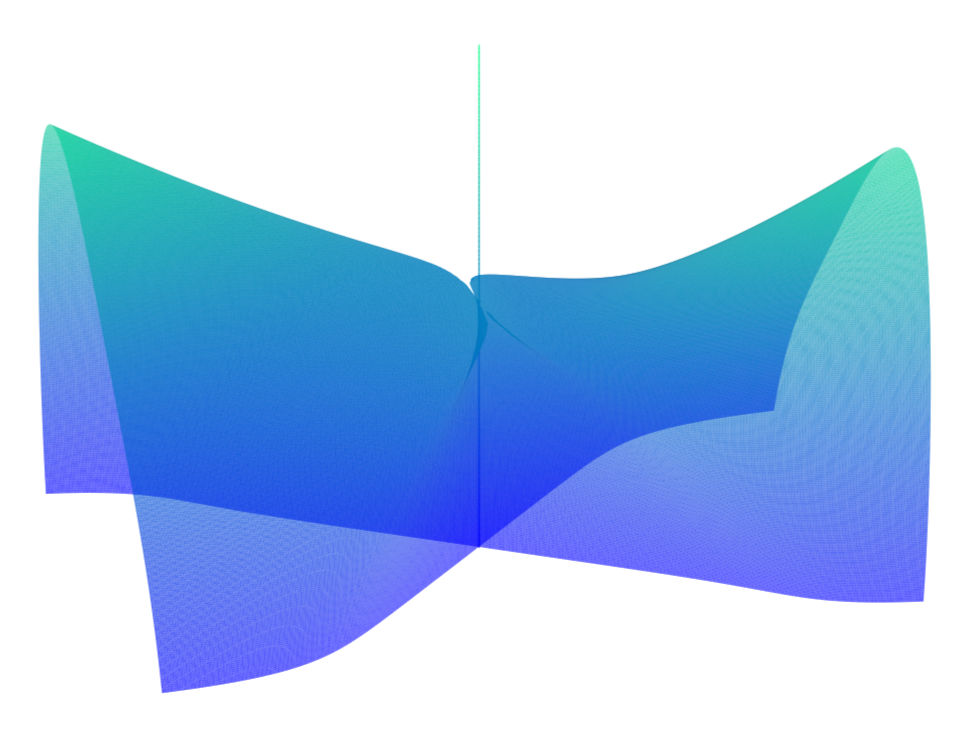}
        \caption{}
        \label{fig:singularities_nonqs}
    \end{subfigure}
    \caption{(\subref{fig:singularities_qs}) The real locus of the affine cone of the member of $X_6\subset\mathbb{P}(1,2,3)$ whose coefficients are all equal to $1$, which is quasismooth and hence admits only quotient singularities. (\subref{fig:singularities_nonqs}) The real locus of the affine cone of the member of $X_{10}\subset\mathbb{P}(1,3,4)$ whose coefficients are all equal to $1$, which is nonquasismooth; the hyperquotient singularity is visible as the line passing through the origin.}
    \label{fig:singularities}
\end{figure} 

We now turn to the analysis of terminal singularities. Singularities of a hypersurface $X$ can arise in two ways: from the vanishing of all partial derivatives of its defining equation $f$, or from nontrivial stabilisers of the $\mathbb{C}^*$-action defining the ambient weighted projective
space. When only the latter occurs, we obtain a \textit{quotient singularity}, whose terminality can be checked efficiently. Quasismooth hypersurfaces have only such singularities, as illustrated in Figure~\ref{fig:singularities} (\subref{fig:singularities_qs}). When only the former occurs, we obtain a \textit{hypersurface singularity}. When both occur simultaneously, we obtain a \textit{hyperquotient singularity}, combining a singular hypersurface with a nontrivial quotient structure. This is illustrated in Figure~\ref{fig:singularities} (\subref{fig:singularities_nonqs}), where the singular locus contains a line through the origin.

We restrict our search to \textit{general} hypersurfaces $X_d\subset\mathbb{P}(a_0,\ldots,a_{n+2})$, corresponding to a dense Zariski-open subset of the parameter space of weighted homogeneous equations of degree $d$. This allows us to determine the singularities occurring generically in each family and analyse their terminality uniformly.

If $P$ is a cyclic quotient singularity \citep[\S~4]{YPG}, then it is locally of type $\tfrac{1}{r}(b_1,\ldots,b_n)$, where $r\geq 1$ and $0\leq b_i\leq r-1$.  The Reid--Shepherd-Barron--Tai criterion \citep[\S3.1]{C3F}\citep[\S3.2]{T} states that $P$ is terminal if and only if
\[
\frac{1}{r}\left(\sum_{i=1}^{n}\overline{kb_i}\right)-1>0,
\qquad \forall\,1\leq k\leq r-1,
\]
where $\overline{kb_i}\in\{0,\ldots,r-1\}$ denotes the residue of $kb_i$ modulo $r$. For $r=1$, the quotient is smooth and hence terminal.

If $P$ is a hyperquotient singularity \citep[\S~4]{YPG}, then it is locally of type $\tfrac{1}{r}(b_1,\ldots,b_{n+1};e)$, where $r\geq 1$ and $0\leq b_i,e\leq r-1$. Since testing terminality directly can be computationally expensive, we first apply an approximation to Mori's criterion \citep{M1}, restricting to lattice points inside the unit cube. For $r\geq 2$, the preliminary test requires
\[
\frac{1}{r}\left(\sum_{i=1}^{n+1}\overline{kb_i}\right)
-\min\left\{
\frac{1}{r}\left(\sum_{i=1}^{n+1}m_i\overline{kb_i}\right)
\ \Bigg\lvert\
x_1^{m_1}\cdots x_{n+1}^{m_{n+1}}\in f'
\right\}-1>0,
\quad \forall\,1\leq k\leq r-1,
\]
where $f'$ is the local analytic equation defining the hypersurface near $P$, and $\overline{kb_i}$ denotes the residue of $kb_i$ modulo $r$. When $r=1$, the preliminary test imposes no condition, and terminality must be determined separately. The approximation provides a necessary but not sufficient condition for terminality and is substantially less computationally expensive to evaluate than the full criterion. We therefore use it as a preliminary test, before certifying the terminality of all reported families using the complete procedure developed in \citep{MTFano, Thesis}. Remarkably, this preliminary test proves highly effective in practice, achieving an accuracy of $99.93\%$, as quantified in \S\ref{sec:search_results}.

\section{Search Algorithms}
\label{sec:integer_lattice_search}

\subsection{Problem formulation}

Our goal is to discover new families of Fano $4$-fold hypersurfaces with terminal singularities. As described in \S\ref{sec:background}, candidate families are represented by weight vectors in the $6$-dimensional integer lattice $\mathbb{Z}^6$, with those satisfying the terminality conditions corresponding to \emph{reward points}. The quasismooth classification, as was depicted in Figure~\ref{fig:class_qs}, suggests that these points, although sparse, exhibit spatial clustering, a property we exploit to guide the search. More generally, locating reward points in $\mathbb{Z}^n$ presents two main challenges:
\begin{enumerate}
    \item \textbf{Sparsity:} Reward points occupy only a small
    fraction of the search space.
    \item \textbf{Unknown a priori:} The reward status of a
    point cannot be determined without directly evaluating it.
\end{enumerate}
However, the reward landscape also exhibits an exploitable property:
\begin{enumerate}
    \item \textbf{Spatial clustering:} Reward points tend to
    occur near one another, so finding a reward point suggests
    that neighbouring points are more likely to yield rewards.
\end{enumerate}

Our objective is to discover both a large number of reward points and those that are difficult to reach using conventional search methods. Since sparsity makes exhaustive exploration expensive, we instead exploit spatial clustering to guide the search. The two algorithms presented below use previously discovered rewards to prioritise promising regions of the lattice, but differ in how these priorities are determined.

\subsection{Fixed heuristic}

The algorithm begins with a start point. We proceed by searching its neighbouring points and determining whether they are reward points. We add the neighbouring points to a search queue and assign them a priority value dependent on their proximity to previously found reward points. The function that computes this priority value is fixed, therefore making it a fixed heuristic algorithm. The algorithm resets by restarting the process with a point in the queue with the highest priority. 

\begin{figure}[!h]
    \centering
    \includegraphics[width=0.8\linewidth]{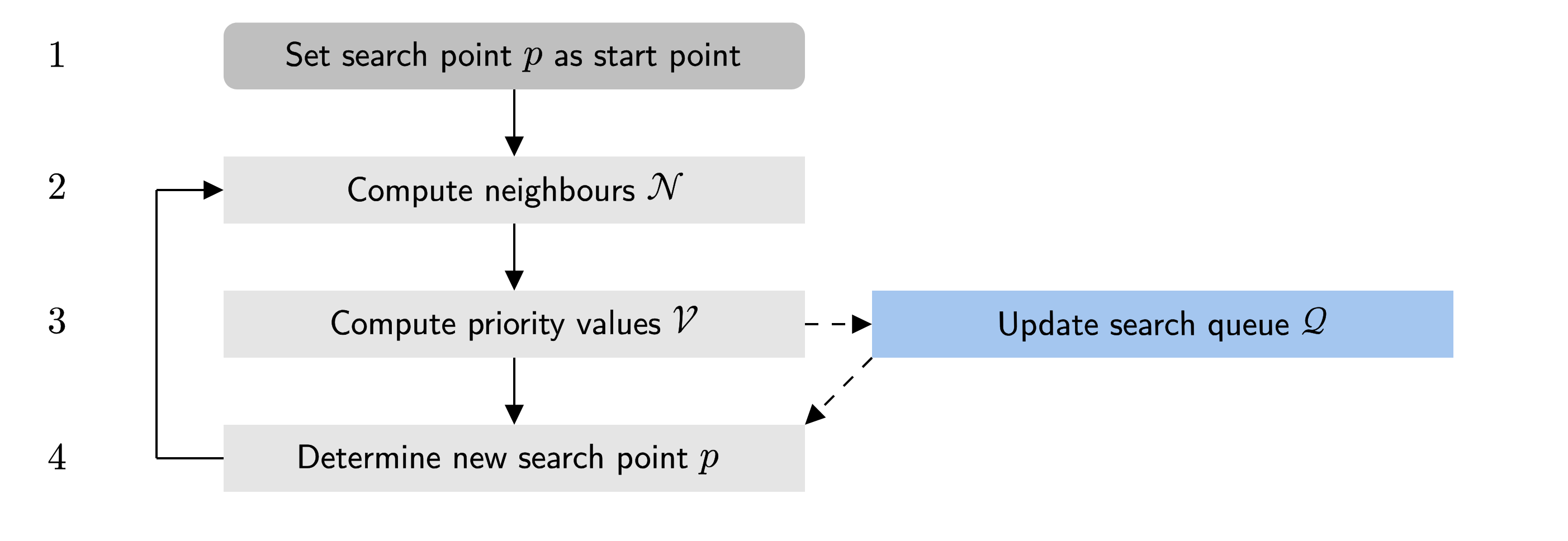}
    \caption[Flowchart of the fixed heuristic search algorithm.]{Flowchart of the fixed heuristic search algorithm.}
    \label{fig:flowchart_fixed}
\end{figure}

 The algorithm depicted in Figure~\ref{fig:flowchart_fixed} is performed as follows.
\begin{enumerate}[left=0pt]
    \item 
    \label{enum:fixed1}
    \begin{enumerate}
        \item Pick a start point $p\in\mathbb{Z}^n$, and set $v(p)=1$, where $v$ is the priority function defined in (\ref{enum:fixed3}).
        \item Set the step count $s=0$ and fix the maximum step count $s_{\text{max}}\in\mathbb{N}$. 
        \item Initialise the search queue $\mathcal{Q}$ as an empty heap.
    \end{enumerate}
    \item 
    \label{enum:fixed2}
    Increment the step count $s$ by $1$. Identify all neighbouring points $\mathcal{N}$ of $p$, defined as the set of points exactly distance $1$ away under the $L^1$ norm, $\mathcal{N}:=\{n\in \mathbb{Z}^n\mid \|n-p\|_1=1\}$, in other words, all points differing by $\pm 1$ from $p$ in one coordinate.  
    \item 
    \label{enum:fixed3}
    For each $n\in\mathcal{N}$ determine its priority value
    \[
    v(n):=
    \begin{cases}
        1, & \text{if $n$ is a reward point,} \\
        \tfrac{1}{2}v(p), & \text{otherwise}, 
    \end{cases}
    \]
    and add $(n,v(n))$ to the search queue $\mathcal{Q}$ if $n$ has never been added before. 
    \item 
    \label{enum:fixed4}
    Determine a point $p'$ such that $(p',v(p'))\in\mathcal{Q}$ has the largest value $v(p')$ in the heap. That is, take the first point of the heap ordered by priority values $v$. Set $p=p'$ and remove $(p',v(p'))$ from $\mathcal{Q}$. Return to (\ref{enum:fixed2}) if $s<s_{\text{max}}$, otherwise terminate the algorithm. 
\end{enumerate}

We observe in \S\ref{sec:results} that when the algorithm is applied to finding terminal Fano hypersurfaces, it is effective at finding many new examples in reward dense regions. We build on the ideas of the fixed heuristic search to design a dynamic heuristic search in \S\ref{sec:algorithm_dynamic} that can find reward points in lower density areas. 

\subsection{Dynamic heuristic (deep reinforcement learning)}
\label{sec:algorithm_dynamic}

The algorithm begins with a chosen start point. We compute its neighbours and, for each, determine the priority values assigned to them by a neural network function. We add these to a search queue ordered by priority values. Next, we assign rewards dependant on whether the neighbours added were reward points or not, and use these to update the neural network using temporal difference learning. The process is then repeated by searching a point with the highest priority in the search queue. 

\begin{figure}[!h]
    \centering
    \includegraphics[width=0.8\linewidth]{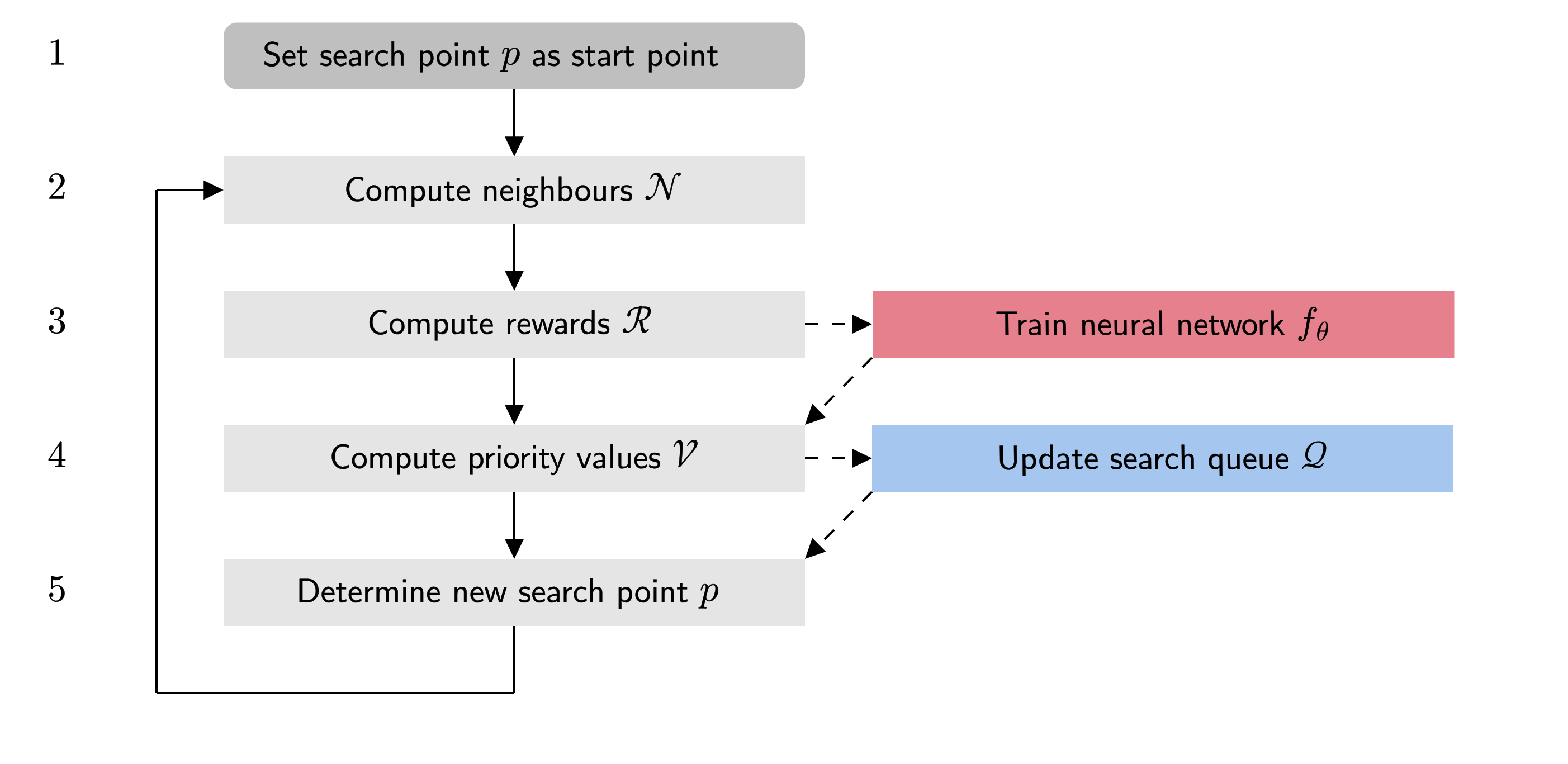}
    \caption[Flowchart of the dynamic heuristic search algorithm.]{Flowchart of the dynamic heuristic search algorithm.}
    \label{fig:flowchart_dynamic}
\end{figure}

The algorithm depicted in Figure \ref{fig:flowchart_dynamic} is performed as follows.
\begin{enumerate}[left=0pt]
    \item 
    \begin{enumerate}
        \item 
        \label{enum:dynamic1}
        Pick a start point $p\in\mathbb{Z}^n$.
        \item Set the step count $s=0$, and fix the maximum step count $s_\text{max}\in\mathbb{N}$.  
        \item Initialise a search queue $\mathcal{Q}$ as an empty heap.
        \item  Create an MLP neural network $f_\theta\colon\mathbb{Z}^n\rightarrow\mathbb{R}$ with initial parameter $\theta$, this will be our dynamic heuristic that determines priority in the search queue $\mathcal{Q}$. Fix the temporal difference discount factor $\gamma\in(0,1)$, which affects how we update $f_\theta$ via temporal difference learning.  
        \item Fix a standard deviation $\sigma\in\mathbb{R}_{\geq 0}$, this determines the stochastic component added to the priority value and thereby controls exploration.
        \item Fix $r_\text{reward}\in\mathbb{N}$, the value given for finding a reward point. Set $s_\text{reward}=0$, the number of steps since a reward was last found.  
    \end{enumerate}
    \item 
    \label{enum:dynamic2}
    Increment the step count $s$ and steps since terminal $s_\text{reward}$ by $1$. Identify all neighbouring points $\mathcal{N}$ of $p$, defined as the set of points exactly distance $1$ away under the $L^1$ norm, $\mathcal{N}:=\{n\in\mathbb{Z}^n\mid \|n-p\|_1=1\}$, in other words, all points differing by $\pm 1$ from $p$ in one coordinate. 
    \item 
    \label{enum:dynamic3}
    Determine whether any $n\in\mathcal{N}$ are reward points, and if so, reset $s_\text{reward}=0$. Compute their reward values
    \[
    r(n)=
    \begin{cases}
        r_\text{reward}, &\text{if $n$ is a reward point},\\
        -\sqrt{s_\text{reward}}, &\text{otherwise}. 
    \end{cases}
    \]
    Consider the set of tuples $(p, n, r(n))$ for each $n \in \mathcal{N}$. This data is used to train the network via temporal difference (TD) learning \citep[\S6]{Sutton2018}. To improve training stability, we fix a copy of the current network parameters, denoting them $\theta'$, which remain frozen during this update step. For each $(p, n, r(n))$ for $n\in \mathcal{N}$, we compute their TD targets
    \[
    t(n) = r(n) + \gamma f_{\theta'}(n),
    \]
    where $\gamma \in (0,1)$ is the discount factor controlling the trade off between short and long term rewards. Values of $\gamma$ close to $0$ produce greedy, short term behaviour whilst values close to $1$ encourage more long term behaviour. We then compute the TD error, measuring the discrepancy between the estimated value of $p$ and the TD target,
    \[
    \delta(\theta, p, n) = f_\theta(p) - t(n).
    \]
    Note that $\theta$ in $f_\theta(p)$ is updated during optimisation, whilst $t(n)$ is held fixed via $\theta'$. Minimising $\delta(\theta, p, n)$ constitutes bootstrapping: future value estimates are refined using past ones. Concretely, we minimise the normalised mean squared error (MSE) loss
    \[
    L(\theta) = \tfrac{1}{2\lvert\mathcal{N} \rvert}\sum_{n \in \mathcal{N}} \delta(\theta,p,n)^2,
    \]
    using a gradient based optimiser such as Adam.
    \item 
    \label{enum:dynamic4}
    For each $n\in\mathcal{N}$, compute their priority values $v(n)=f_\theta(n)+\varepsilon$, where $\varepsilon$ is sampled from $\mathcal{N}(0,\sigma^2)$, a normal distribution with mean $0$ and variance $\sigma^2$. The stochastic component, $\varepsilon$, improves exploration. Add $(n,v(n))$ to the search queue $\mathcal{Q}$ if it has not previously been searched before.
    \item 
    \label{enum:dynamic5}
    Let $p'$ be a point in the search queue $\mathcal{Q}$ such that $v(p')$ is the largest value in the heap. That is, take the first point of the heap ordered by priority values $v$. Set the new search point $p=p'$, and return to (\ref{enum:dynamic2}) if $s<s_{\text{max}}$, otherwise terminate the algorithm. 
\end{enumerate}

Since the dynamic heuristic search is nondeterministic, rerunning the algorithm can uncover new rewards within the same fixed number of steps. An ablation study in \S\ref{sec:ablation} confirms that the neural network is responsible for the ability to find hard-to-reach examples, with the stochastic component alone insufficient to reach them.

\section{Computational results}
\label{sec:results}

\subsection{Experimental setup}
\label{sec:experimental_setup}

We apply the fixed and dynamic heuristics of \S\ref{sec:integer_lattice_search} to discover nonquasismooth families of terminal Fano $4$-fold hypersurfaces of index~$1$. Our primary interest is in discovering new families with degree $d>200$, beyond the range of the exhaustive classification discussed in \S\ref{sec:computational_obstruction}. Our goal is twofold: to identify as many new examples as possible and to uncover hard-to-reach ones. We show that the fixed heuristic search is particularly successful in the former, whilst the dynamic heuristic search achieves both. 

To initialise the searches, we use the $11{,}617$ quasismooth terminal families classified by Brown and Kasprzyk~\citep{BK}. Specifically, the first $11{,}617$ points processed by each algorithm are these known terminal examples, providing an initial collection of rewards from which the searches can proceed. Thereafter, points are selected according to the respective heuristics. Each search is run for $10{,}000{,}000$ steps, including the initialisation. The fixed heuristic is deterministic and is run once, while the dynamic heuristic is run five times independently. The hyperparameters used in the dynamic search are given in Table~\ref{tab:hyperparameters}.

\begin{table}[!h]
    \centering
    \small
    \setlength{\arrayrulewidth}{3pt}
    \arrayrulecolor{white}
    \begin{tabular}{@{}ll@{\hspace{2em}}ll@{}}
        \rowcolor{gray!20}
        \textbf{Hyperparameter} & \textbf{Value}
        & \textbf{Hyperparameter} & \textbf{Value} \\ \hline

        MLP hidden layers & \cellcolor{gray!10}$(40,)$
        & Learning rate & \cellcolor{gray!10}$0.001$ \\

        Activation & \cellcolor{gray!10}LeakyReLU
        & TD discount $\gamma$ & \cellcolor{gray!10}$0.2$ \\

        LeakyReLU slope & \cellcolor{gray!10}$0.01$
        & Standard deviation $\sigma$ & \cellcolor{gray!10}$2$ \\

        Optimiser & \cellcolor{gray!10}Adam
        & Search reward $r_{\mathrm{reward}}$
        & \cellcolor{gray!10}$1$ \\
    \end{tabular}
    \caption{Hyperparameters used in the dynamic heuristic
    (deep reinforcement learning) search.}
    \label{tab:hyperparameters}
\end{table}

The algorithms were implemented in Python, using PyTorch 2.8.0~\citep{pytorch} for neural network training. Experiments were conducted on a laptop with an Intel Core i7-1185G7 processor and 16 GB RAM, using CPU only. For $10{,}000{,}000$ steps, the fixed and dynamic heuristics took approximately $1.5$ and $3$ hours, respectively, with subsequent dynamic runs benefiting from cached evaluations.

\subsection{Distance based analysis of search reach}
\label{sec:distance_analysis}

Let $\mathcal{Q}$ denote the set of $11{,}617$ quasismooth terminal families classified by Brown and Kasprzyk. Let $\mathcal{F}$ denote the set of distinct nonquasismooth terminal families discovered by the fixed heuristic, and $\mathcal{D}_i$ the corresponding set from the $i$-th dynamic run, for $1\leq i\leq 5$. To quantify how isolated the discovered families are, we measure the $L^1$-distance from each point to its nearest distinct point among the corresponding search results and the known quasismooth families. For each $\mathcal{S}\in\{\mathcal{F},\mathcal{D}_1,\ldots, \mathcal{D}_5\}$ and $p\in\mathcal{S}$, we define
\[
    D_{\mathcal{S}}(p):=\min\{\|p-q\|_1 \mid q\in(\mathcal{S}\sqcup\mathcal{Q})\setminus\{p\}\}.
\]
Thus, $D_{\mathcal{S}}(p)$ measures the $L^1$-distance from $p$ to its nearest distinct point in the corresponding search results or the quasismooth classification. From these distances we derive a lower bound on the number of steps required to reach a point from its nearest neighbour. This allows us to show that points found by the dynamic search would be computationally expensive to find using the fixed search alone. 

Let $q$ be a nearest neighbour of $p$, so that $D_{\mathcal{S}}(p)=\|p-q\|_1$. We estimate the number of steps required for the fixed heuristic to discover $p$ starting from $q$. Since the fixed heuristic halves the priority at each non-reward point, exploration proceeds in increasing $L^1$-distance from $q$ unless another reward changes the priorities. In the favourable case where exploration is not redirected, the search must expand all points within distance $D_{\mathcal{S}}(p)-2$ of $q$, together with at least one point at distance $D_{\mathcal{S}}(p)-1$, before discovering $p$. Conversely, expanding all points within distance $D_{\mathcal{S}}(p)-1$ guarantees its discovery. Writing $\Lambda:=\{(a_1,\ldots,a_6)\in\mathbb{Z}^6: 1\leq a_1\leq\cdots\leq a_6\}$, we obtain
\[
\begin{aligned}
s_L(p):=\#\bigl(B_1(q,D_{\mathcal{S}}(p)-2)\cap\Lambda\bigr)+1,\quad s_U(p):=\#\bigl(B_1(q,D_{\mathcal{S}}(p)-1)\cap\Lambda\bigr),
\end{aligned}
\]
where $B_1(q,r)$ denotes the closed $L^1$-ball of radius $r$ centred at $q$.

The lower bound $s_L(p)$ is deliberately conservative, as it assumes the most favourable exploration towards $p$. In practice, competing rewards may redirect the fixed heuristic, increasing the number of steps required. Consequently, while $s_L(p)$ remains a lower bound under the fixed heuristic's priority ordering, $s_U(p)$ is an upper bound only when exploration from $q$ is not redirected.

To estimate the likelihood of attaining the lower bound, we consider an idealised probabilistic model in which the points at distance $D_{\mathcal{S}}(p)-1$ from $q$ are explored in uniformly random order, with exactly one leading to $p$. The probability of discovering $p$ within $s_L(p)\leq s\leq s_U(p)$ steps is then
\[
P(s)=\frac{s-s_L(p)+1}{s_U(p)-s_L(p)+1}.
\]
In particular, $P(s_L(p))=1/(s_U(p)-s_L(p)+1)$. As the number of points in this distance shell grows with $D_{\mathcal S}(p)$, attaining the lower bound becomes increasingly unlikely. This model quantifies the conservativeness of the lower bound, rather than predicting the actual discovery time of the fixed heuristic.

The exact lower bound depends on how the ball $B_1(q,r-2)$ intersects the boundary of $\Lambda$. Since almost all examples have first weight $1$, the dominant boundary effect is the constraint $a_1\geq1$. Further truncation by the remaining boundary faces $a_i=a_{i+1}$ was found to be negligible in practice. Assuming $q_1=1$ and neglecting these additional boundary effects, we obtain the approximation
\[
\widetilde{s}_L(p)
:=\#\bigl(B_1(0,D_{\mathcal S}(p)-2)\cap
\{a\in\mathbb{Z}^6:a_1\geq0\}\bigr)+1.
\]
Table~\ref{tab:distance_cost} illustrates the rapid growth of the estimated search cost with nearest-neighbour distance. In particular, at distance $22$, the conservative estimate exceeds $3{,}800{,}000$ steps, almost $40\%$ of the $10{,}000{,}000$ step search budget.

\begin{table}[!h]
    \centering
    \small
    \setlength{\arrayrulewidth}{3pt}
    \arrayrulecolor{white}
    \begin{tabular}{@{}llllllllll@{}}
        \rowcolor{gray!20}$\boldsymbol{D_{\mathcal{S}}(p)}$
        & $1$
        & $5$
        & $10$
        & $\cdots$
        & $17$
        & $18$
        & $19$
        & $21$
        & $22$ \\ \hline
        \rowcolor{gray!10}$\boldsymbol{\widetilde{s}_L(p)}$
        & $1$
        & $305$
        & $26{,}578$
        & $\cdots$
        & $774{,}913$
        & $1{,}110{,}050$
        & $1{,}558{,}477$
        & $2{,}915{,}761$
        & $3{,}898{,}490$ \\
    \end{tabular}
    \caption{Estimated lower-bound search cost as a function
    of the nearest-neighbour distance.}
    \label{tab:distance_cost}
\end{table}

For example, $X_{1020}\subset\mathbb{P}(1,15,32,139,340,494) \in\mathcal{D}_1\setminus\mathcal{F}$ has nearest-neighbour distance $D_{\mathcal{D}_1}(p)=17$, giving $\widetilde{s}_L(p)=774{,}913$. Its nearest neighbour is $X_{1011}\subset\mathbb{P}(1,10,31,143,337,490)\in \mathcal{Q}$. Although a constant-priority search, where we set $v(n)=1$, finds this example in $1{,}041{,}501$ steps, the fixed heuristic fails to reach it within $10{,}000{,}000$ steps, even when initialised from its nearest neighbour alone. A more isolated example is $X_{1220}\subset\mathbb{P}(1,9,32,173,396,610)\in \mathcal{D}_2\setminus \mathcal{F}$, whose nearest neighbour is $X_{1204}\subset\mathbb{P}(1,7,24,172,399,602) \in\mathcal{Q}$ at distance $D_{\mathcal{D}_2}(p)=22$. Here, $\widetilde{s}_L(p)=3{,}898{,}490$, almost $40\%$ of the entire $10{,}000{,}000$ step budget. These examples illustrate the difficulty of reaching isolated families using the fixed heuristic.

\subsection{Heuristic search results}
\label{sec:search_results}

\begin{figure}[!h]
    \centering
    \includegraphics[width=0.8\linewidth]{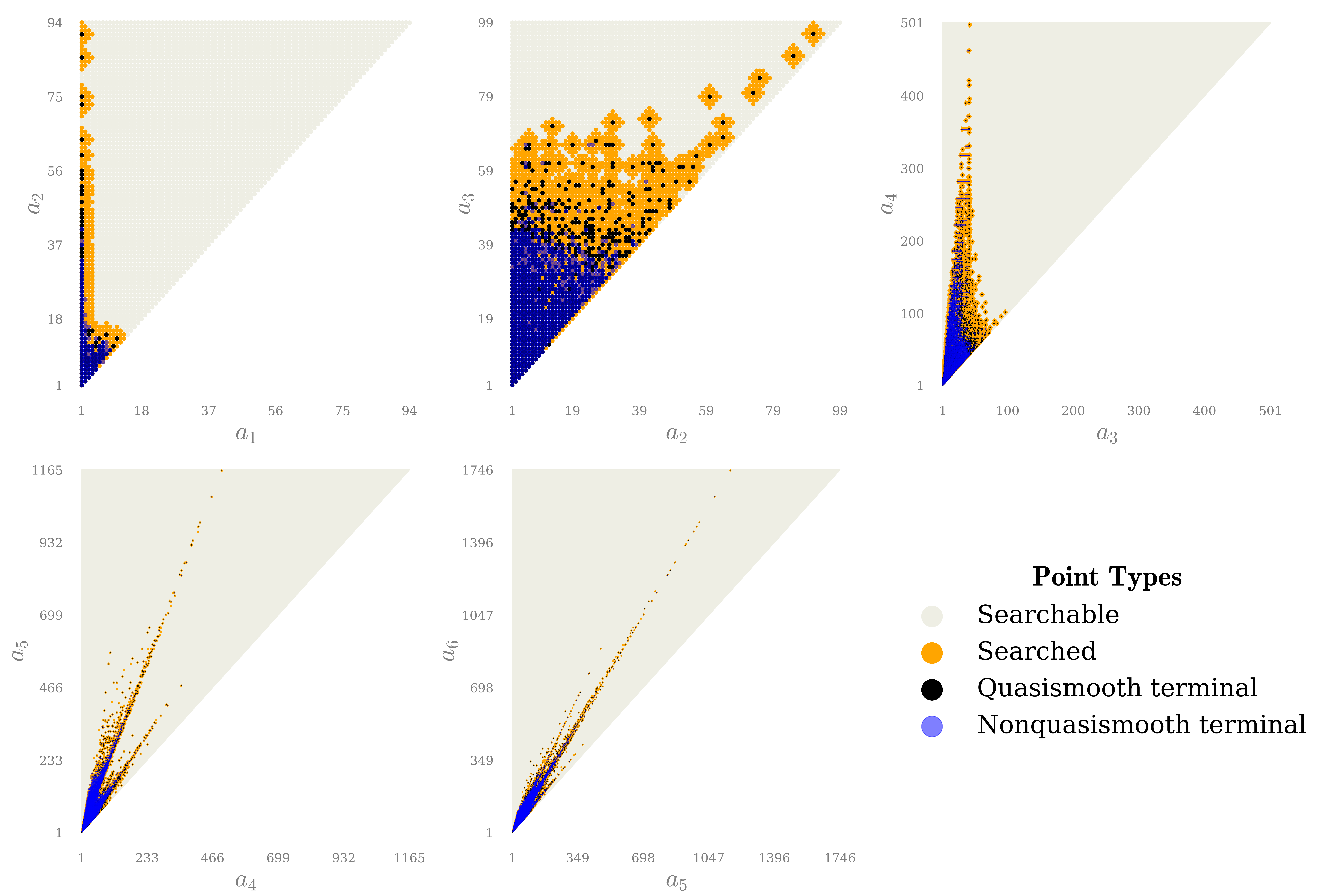}
    \caption{Fixed heuristic search, corresponding to $\mathcal{F}$. Each panel shows a projection of the search from $\mathbb{Z}^6$ onto $\mathbb{Z}^2$ via $(a_1,\ldots,a_6)\mapsto(a_i,a_{i+1})$, illustrating the distribution of discovered families across consecutive weight coordinates.}
    \label{fig:search_fixed}
\end{figure}

\begin{figure}[!h]
    \centering
    \includegraphics[width=0.8\linewidth]{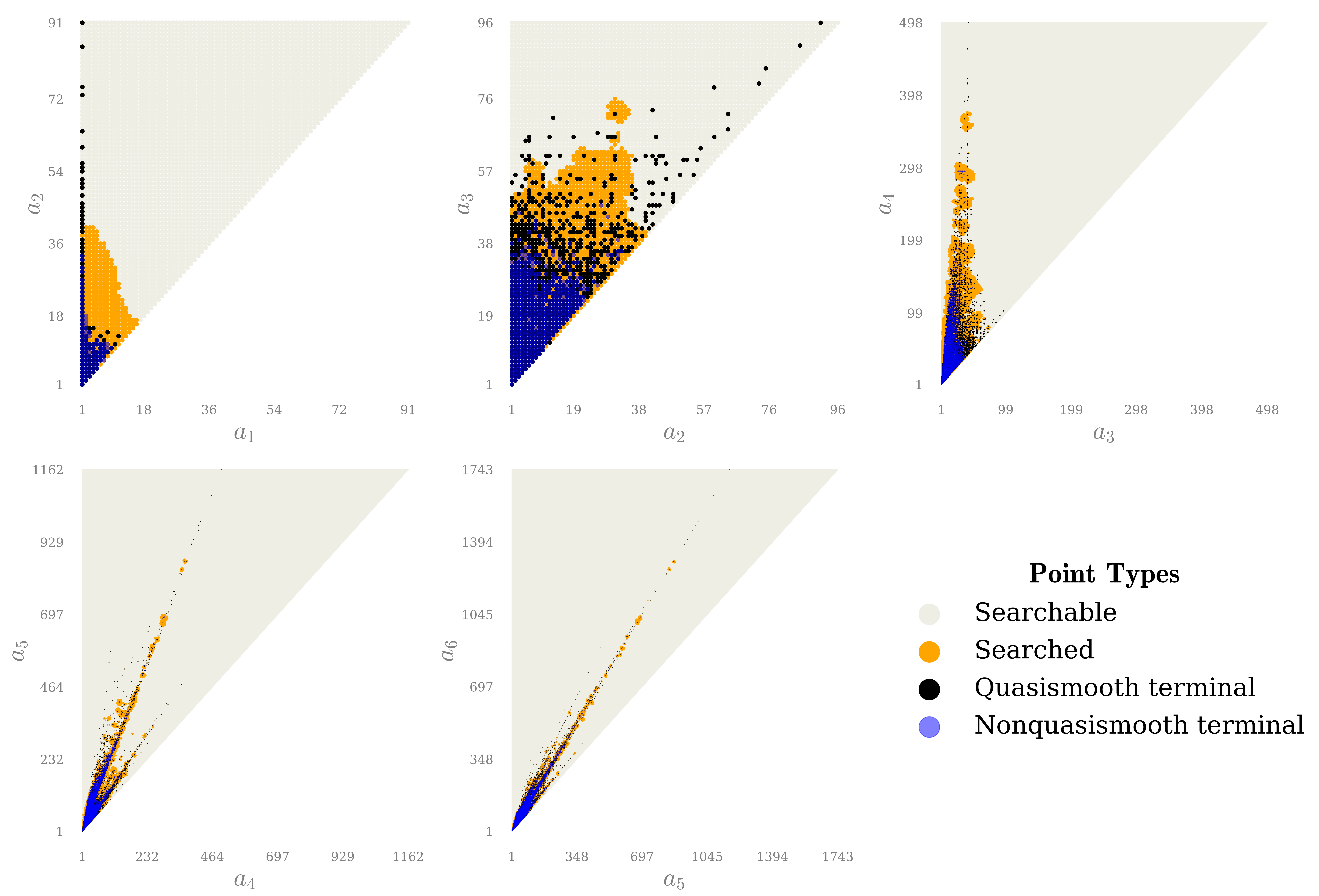}
    \caption{Dynamic heuristic (deep reinforcement learning) search, corresponding to the first run $\mathcal{D}_1$. Each panel shows a projection of the search from $\mathbb{Z}^6$ onto $\mathbb{Z}^2$ via $(a_1,\ldots,a_6)\mapsto(a_i,a_{i+1})$, illustrating the distribution of discovered families across consecutive weight coordinates.}
    \label{fig:search_dynamic}
\end{figure}

As illustrated in Figures~\ref{fig:search_fixed} and \ref{fig:search_dynamic}, the fixed heuristic and the first dynamic run exhibit distinct exploration patterns. Table~\ref{tab:search_results} compares their results and summarises the combined discoveries $\mathcal{U}=\mathcal{F}\cup\bigcup_{i=1}^{5}\mathcal{D}_i$. To quantify the exploration of more isolated regions, we call a family \emph{hard to reach} if $d>200$ and $D_{\mathcal S}(p)\geq10$, where $D_{\mathcal S}(p)$ is the distance statistic defined in \S\ref{sec:distance_analysis}.

\begin{table}[!h]
    \centering
    \small
    \setlength{\tabcolsep}{2pt}
    \renewcommand{\arraystretch}{1.2}
    \setlength{\arrayrulewidth}{3pt}
    \arrayrulecolor{white}

    \begin{tabular}{@{}
        >{\centering\arraybackslash}m{0.07\linewidth}
        >{\centering\arraybackslash}m{0.08\linewidth}
        >{\centering\arraybackslash}m{0.08\linewidth}
        >{\centering\arraybackslash}m{0.30\linewidth}
        >{\centering\arraybackslash}m{0.12\linewidth}
        >{\centering\arraybackslash}m{0.26\linewidth}
    @{}}
        \rowcolor{gray!40}
        $\boldsymbol{\mathcal S}$
        & $\boldsymbol{|\mathcal S|}$
        & $\boldsymbol{d>200}$
        & $\boldsymbol{\max d}$
        & $\boldsymbol{D_{\mathcal S}(p)\geq10}$
        & $\boldsymbol{\max D_{\mathcal S}(p)}$ \\ \hline

        \rowcolor{gray!10}
        $\mathcal F$
        & $113{,}922$
        & $41{,}333$
        & {\scriptsize $X_{3234}\subset\mathbb P(1,37,40,462,1078,1617)$}
        & $0$
        & \shortstack[c]{\\[2pt]
            $6$\\[3pt]
            {\scriptsize $X_{367}\subset\mathbb P(1,2,21,73,98,173)$}\\[1pt]
            {\scriptsize $X_{364}\subset\mathbb P(1,1,28,48,105,182)$} \\[0pt]
          } \\ \hline

        \rowcolor{gray!20}
        $\mathcal D_1$
        & $85{,}214$
        & $18{,}693$
        & {\scriptsize $X_{2070}\subset\mathbb P(1,19,33,293,690,1035)$}
        & $91$
        & \shortstack[c]{\\[2pt]
            $17$\\[3pt]
            {\scriptsize $X_{1020}\subset\mathbb P(1,15,32,139,340,494)$} \\[0pt]
          } \\

        \rowcolor{gray!20}
        $\mathcal D_2$
        & $81{,}053$
        & $16{,}049$
        & {\scriptsize $X_{3234}\subset\mathbb P(1,37,40,462,1078,1617)$}
        & $82$
        & \shortstack[c]{\\[2pt]
            $22$\\[3pt]
            {\scriptsize $X_{1220}\subset\mathbb P(1,9,32,173,396,610)$} \\[0pt]
          } \\

        \rowcolor{gray!20}
        $\mathcal D_3$
        & $84{,}102$
        & $17{,}750$
        & {\scriptsize $X_{2604}\subset\mathbb P(1,23,39,372,868,1302)$}
        & $69$
        & \shortstack[c]{\\[2pt]
            $18$\\[3pt]
            {\scriptsize $X_{868}\subset\mathbb P(1,6,25,124,279,434)$} \\[0pt]
          } \\

        \rowcolor{gray!20}
        $\mathcal D_4$
        & $86{,}676$
        & $18{,}469$
        & {\scriptsize $X_{2604}\subset\mathbb P(1,23,39,372,868,1302)$}
        & $71$
        & \shortstack[c]{\\[2pt]
            $22$\\[3pt]
            {\scriptsize $X_{906}\subset\mathbb P(1,21,31,99,302,453)$} \\[0pt]
          } \\

        \rowcolor{gray!20}
        $\mathcal D_5$
        & $72{,}068$
        & $12{,}201$
        & {\scriptsize $X_{2772}\subset\mathbb P(1,25,41,396,924,1386)$}
        & $82$
        & \shortstack[c]{\\[2pt]
            $20$\\[3pt]
            {\scriptsize $X_{1960}\subset\mathbb P(1,22,37,280,641,980)$}\\[3pt]
            {\scriptsize $X_{1122}\subset\mathbb P(1,26,36,137,362,561)$} \\[0pt]
          } \\

        \rowcolor{gray!30}
        \shortstack[c]{
            {\scriptsize\bfseries Mean}\\[-2pt]
            {\scriptsize $\boldsymbol{\pm}$ std}
        }
        & {\scriptsize $81{,}823\pm5{,}831$}
        & {\scriptsize $16{,}632\pm2{,}685$}
        & $2{,}657\pm417$
        & $79\pm9$
        & $20\pm2$ \\ \hline

        \rowcolor{gray!40}
        $\boldsymbol{\mathcal U}$
        & $\boldsymbol{121{,}980}$
        & $\boldsymbol{46{,}992}$
        & {\scriptsize $\mathbf{ X_{3234}\subset\mathbb P(1,37,40,462,1078,1617)}$}
        & ---
        & --- \\
    \end{tabular}

    \caption{Results of the fixed heuristic and five independent dynamic runs after $10{,}000{,}000$ steps. Counts refer to distinct nonquasismooth terminal families. The columns report total discoveries, those with $d>200$, and maximum degree. Hard-to-reach families are defined as those satisfying $d>200$ and $D_{\mathcal S}(p)\geq10$; the final column reports the maximum nearest-neighbour distance among families with $d>200$. Examples attaining the maxima are displayed using $X_d\subset\mathbb P(a_1,\ldots,a_6)$. The mean $\pm$ sample standard deviation is calculated across dynamic runs.}
    \label{tab:search_results}
\end{table}

The fixed heuristic discovers $113{,}922$ nonquasismooth terminal families, compared with $81{,}823\pm5{,}831$ for the dynamic heuristic across five independent runs. The preliminary terminality criterion also proves highly effective: of the $122{,}063$ distinct candidates passing the approximate test, $121{,}980$ were subsequently certified terminal, with only $83$ failing full certification. This corresponds to a certification rate of $99.93\%$.

\begin{figure}[!h]
    \centering
    \includegraphics[width=0.81\linewidth]{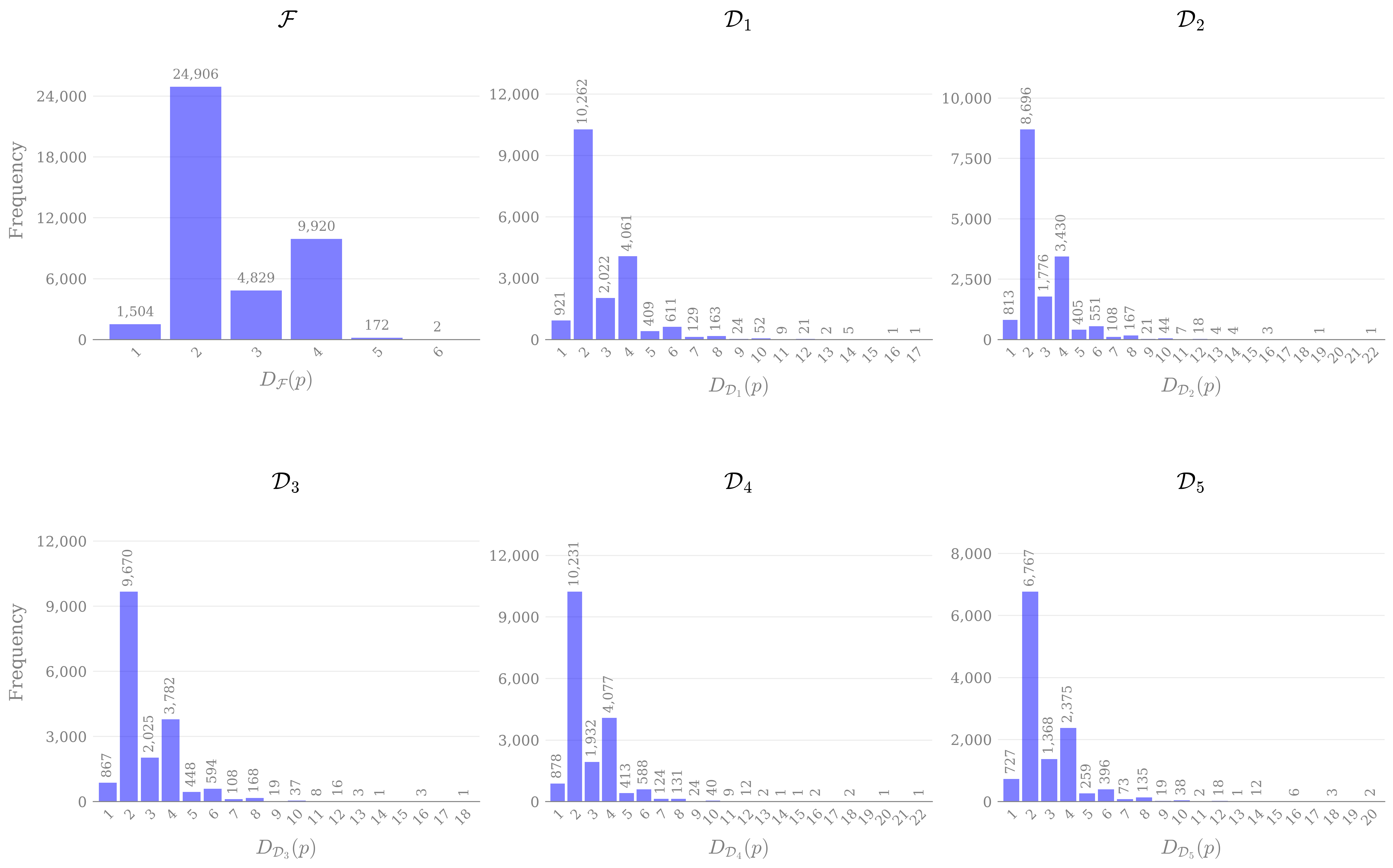}
    \caption{Histograms for the fixed heuristic search $\mathcal F$ and five independent dynamic runs $\mathcal D_i$, showing the distribution of $D_{\mathcal S}(p)$, the $L^1$-distance from each discovered nonquasismooth terminal hypersurface $p$ of degree $d>200$ to its nearest distinct terminal neighbour.  Compared with the fixed heuristic search, the dynamic runs consistently discover families at substantially greater nearest-neighbour distances, including hard-to-reach families satisfying $D_{\mathcal S}(p)\geq10$.}
    \label{fig:histogram}
\end{figure}

Despite its lower overall yield, the dynamic heuristic reaches more isolated regions of weight space. Figure~\ref{fig:histogram} compares the nearest-neighbour distance distributions for families with $d>200$, showing that the dynamic runs consistently reach larger distances. On average, they discover $79\pm9$ hard-to-reach families per run, including examples missed by the fixed search. For the most isolated families, the conservative model developed in \S\ref{sec:distance_analysis} gives estimated lower bounds approaching $40\%$ of the fixed search's step budget to discover a single such family, as detailed in Table~\ref{tab:distance_cost}. This advantage persists across all five runs, and the ablation study in \S\ref{sec:ablation} attributes it to the learned heuristic rather than stochastic exploration alone.

\subsection{Ablation}
\label{sec:ablation}

To determine whether the dynamic heuristic's ability to reach isolated examples arises from the neural network or stochastic exploration alone, we conduct an ablation study. We remove the neural network while retaining the stochastic component, replacing the learned priority function with
\[
    v(n):=
    \begin{cases}
        r_\text{reward}+\varepsilon, & \text{if $n$ is a reward point},\\
        -\sqrt{s_\text{reward}}+\varepsilon, & \text{otherwise},
    \end{cases}
\]
where $\varepsilon\sim\mathcal{N}(0,\sigma^2)$. All other applicable hyperparameters remain as in Table~\ref{tab:hyperparameters}. Let $\mathcal{A}_i$ denote the set of distinct nonquasismooth terminal families discovered by the $i$-th ablation run, for $1\leq i\leq 5$.

\begin{figure}[!htbp]
    \centering
    \begin{minipage}[c]{0.55\linewidth}
        \centering
        \small
        \arrayrulecolor{white}
        \setlength{\arrayrulewidth}{3pt}
        \begin{tabular}{@{}cccc@{}}
            \rowcolor{gray!30}
            $\boldsymbol{\mathcal{A}_i}$
            & $\boldsymbol{|\mathcal{A}_i|}$
            & $\boldsymbol{d>200}$
            & $\boldsymbol{D_{\mathcal{A}_i}(p)\geq10}$ \\ \hline
        
            \rowcolor{gray!10}
            $\mathcal{A}_1$ & $114{,}545$ & $41,076$ & $1$ \\
        
            \rowcolor{gray!10}
            $\mathcal{A}_2$ & $114{,}685$ & $41,117$ & $1$ \\
        
            \rowcolor{gray!10}
            $\mathcal{A}_3$ & $114{,}174$ & $40,684$ & $1$ \\
        
            \rowcolor{gray!10}
            $\mathcal{A}_4$ & $114{,}244$ & $40,795$ & $2$ \\
        
            \rowcolor{gray!10}
            $\mathcal{A}_5$ & $114{,}435$ & $40,995$ & $2$ \\
        
            \rowcolor{gray!20}
            \textbf{Mean $\boldsymbol{\pm}$ std}
            & {\scriptsize $114{,}417\pm211$}
            & {\scriptsize $40,945\pm203$}
            & $1.4\pm0.5$ \\
        \end{tabular}
        \vspace{1.9em}
    \end{minipage}
    \hspace{0.1em}
    \begin{minipage}[c]{0.37\linewidth}
        \centering
        \includegraphics[width=\linewidth]
        {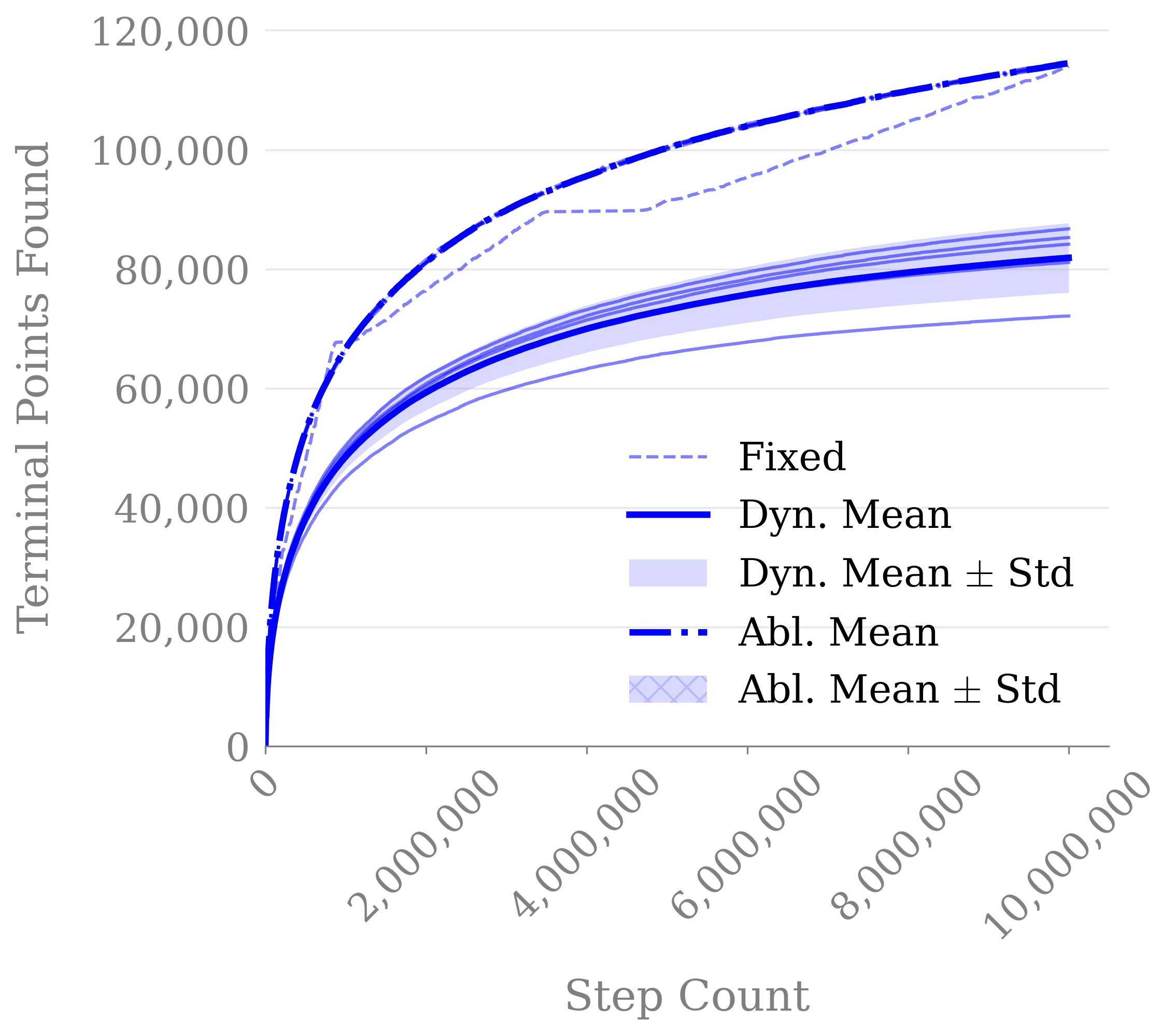}
    \end{minipage}
    \caption{Ablation study over $10{,}000{,}000$ steps. Left: total discoveries for each ablation run. Right: cumulative nonquasismooth terminal families discovered by the fixed, dynamic and ablation searches, with mean $\pm$ standard deviation across five runs for the latter two methods. The ablation curves are nearly indistinguishable.}
    \label{fig:ablation}
\end{figure}

\begin{figure}[!h]
    \centering
    \includegraphics[width=0.81\linewidth]{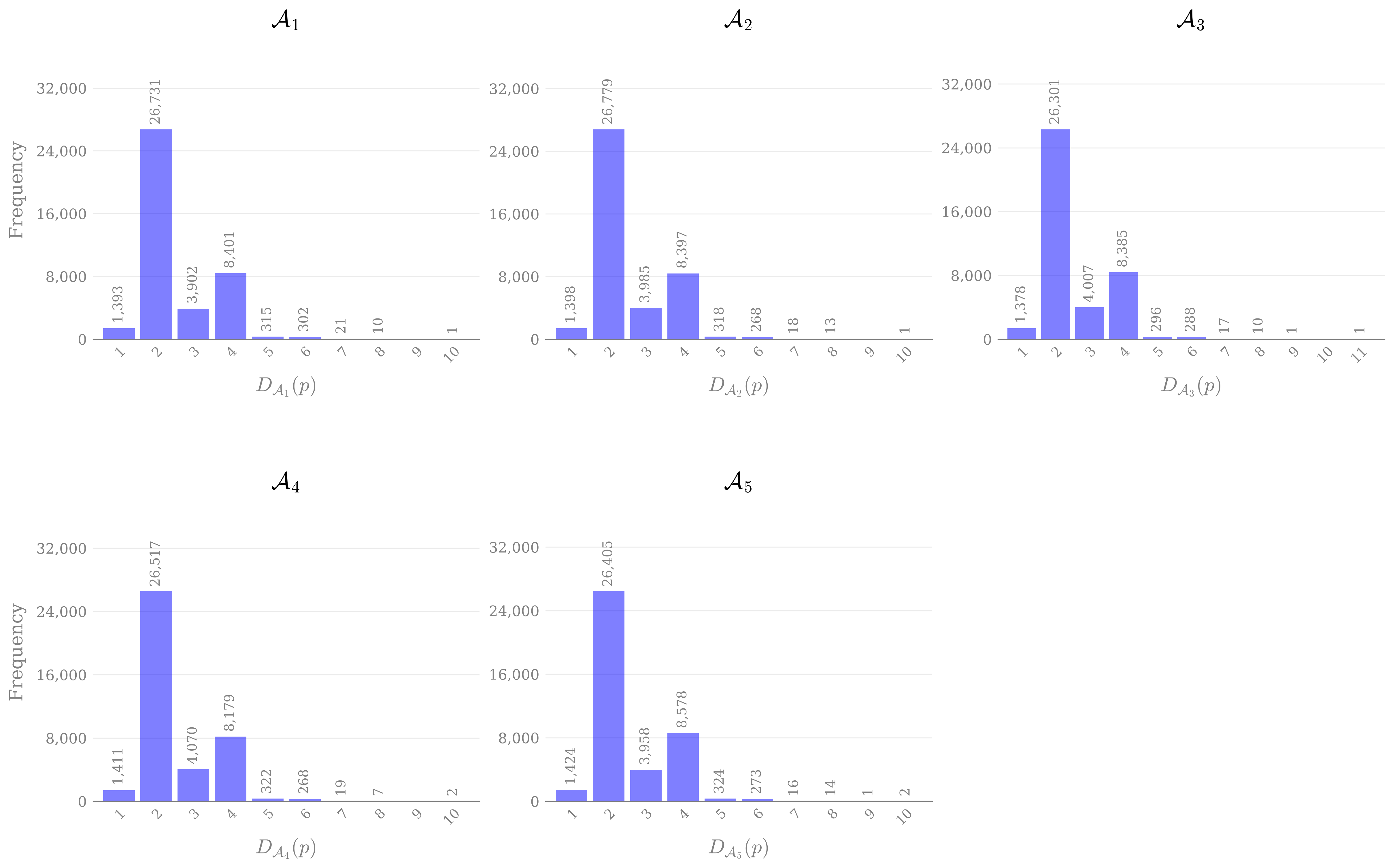}
    \caption{Histograms for the five ablation runs $\mathcal A_i$, showing the distribution of $D_{\mathcal A_i}(p)$, the $L^1$-distance from each discovered nonquasismooth terminal hypersurface $p$ of degree $d> 200$ to its nearest distinct terminal neighbour. Compared with the dynamic heuristic search, the ablation runs find far fewer, if any, families at larger nearest-neighbour distances, particularly those satisfying $D_{\mathcal A_i}(p)\geq 10$.}
    \label{fig:histogram_ablation}
\end{figure}

Across five independent runs, the ablation discovers $114,417\pm211$ terminal families, comparable to the fixed heuristic, although its discovery rate declines more rapidly as seen in Figure~\ref{fig:ablation}. However, the distance distributions in Figure~\ref{fig:histogram_ablation} show that the ablation finds few, if any, hard-to-reach examples, with maximum distances substantially below those of the dynamic heuristic. This provides strong evidence that the neural network, rather than stochastic exploration alone, enables the dynamic heuristic to reach more isolated regions of weight space.

\section*{Acknowledgements}
I am grateful to Gavin Brown, Alexander Kasprzyk, Hefin Lambley and Martin Lotz for valuable feedback during the preparation of this paper. The author was supported by the Warwick Mathematics Institute Centre for Doctoral Training and the UK Engineering and Physical Sciences Research Council (EPSRC) [grant number EP/W523793/1].

{
\bibliographystyle{plainnat}
\bibliography{references}}

\end{document}